\documentclass[12pt, letterpaper]{article}
\usepackage{amsmath}
\usepackage{amsfonts}
\usepackage{amssymb}
\usepackage{amsthm}
\usepackage{breqn}
\usepackage{setspace}
\usepackage{fullpage}
\usepackage{enumerate}
\usepackage{comment}
\usepackage{hyperref}
\usepackage{bbm}
\usepackage{tikz}
\usepackage{mathtools}
\usepackage[utf8]{inputenc}
\usepackage[english]{babel}
\usetikzlibrary{patterns,arrows,decorations.pathreplacing}
\newtheorem{theorem}{Theorem}

\newcommand{\s}{{\rm  s}}

\newcommand{\floor}[1]{\left\lfloor{#1}\right\rfloor}

\title {A note on universal graphs for spanning trees}

\author{Ervin Gy\H{o}ri\thanks{Alfr\'{e}d R\'{e}nyi Institute of Mathematics, Budapest, Hungary.} \and
Binlong Li\thanks{School of Mathematics and Statistics,
Northwestern Polytechnical University, Xi’an, 710072, China.} \thanks{Xi'an-Budapest Joint Research Center for Combinatorics, Northwestern Polytechnical University, Xi'an, 710072, China.}
\and
Nika Salia\thanks{King Fahd University of Petroleum and Minerals, Dhahran, Saudi Arabia.} 
\and 
Casey Tompkins\footnotemark[1]}

\date{}

\begin{document}
\maketitle
\begin{abstract}
Chung and Graham considered the problem of minimizing the number of edges in an $n$-vertex graph containing all $n$-vertex trees as a subgraph.  They showed that such a graph has at least $\frac{1}{2}n \log{n}$ edges.  In this note, we improve this lower estimate to $n \log{n}$. 
\end{abstract}

\section{Introduction}

We say that an $n$-vertex graph $G$ is universal for spanning trees if $G$ contains every $n$-vertex tree as a subgraph.
Chung and Graham~\cite{chung1978graphs} considered the problem of how few edges a graph $G$ can have if it is universal for spanning trees; let us denote this minimum by $\s(n)$. They also considered a related parameter $s^*(n)$ defined to be the minimum number of edges of a graph  (with no condition on the number of vertices) containing every $n$-vertex tree as a subgraph. 
Obviously $\s^*(n) \le \s(n)$.    
A lower bound of $\frac{1}{2}n\log(n)\le {\s}^*(n)$ was given in~\cite{chung1978graphs}, where $\log$ denotes the natural logarithm (see also~\cite{chung1979universal}). Upper bounds for ${\s}^*(n)$ or $\s(n)$ were given in~\cite{chung1978graphs}, \cite{chung1976graphs} and \cite{chung1983universal} yielding

\[ \frac{1}{2}n\log{n} <  {\s}^*(n) \le \s(n) \le \frac{5}{\log{4}}n\log{n} + O(n).\]

In regard to the correct leading coefficient for $\s(n)$, Chung and Graham write in~\cite{chung1983universal}: ``Conceivably, the right value might even be $\frac{1}{2}$ although at present, we certainly do not see how to prove this". 
We show that this is, in fact, not the case by improving the lower bound on $\s(n)$ by a factor of $2$.  
The proof resembles the proof of the lower bound for ${\s}^*(n)$ given in~\cite{chung1978graphs} and~\cite{chung1983universal}.  
However, we keep more careful track of the edges rather than immediately considering the degree sequence.  

\begin{theorem}\label{main}
For every positive integer $n$, 
\[n\log{n} - O(n) \le {\s}^*(n).\]
\end{theorem}

\section{Proof of Theorem \ref{main}}
Let $G$ be a graph which contains all $n$-vertex trees.  Observe that since $G$ contains all $n$-vertex trees as a subgraph, it must also contain all forests with at most $n$-vertices. 

For every integer $i$, such that $1\leq i\leq \frac{n}{2}$, let $F_i$ be a star forest with $i$ connected components each isomorphic to a star with $\floor{\frac{n}{i}}$ vertices.

The proof proceeds as follows: we will sequentially consider forests $F_i$, and at each step, we will color $\floor{\frac{n}{i}}-1$ edges red that have not been colored previously.

Let $F_1'$ be a subgraph of $G$ isomorphic to $F_1$. Let $v_1$ be the central vertex of $F'_1$.  
We color all edges of $F_1'$ red.
Let $F_2'$ be a subgraph of $G$ isomorphic to $F_2$. 
Thus $G$ contains two vertex disjoint stars of $\floor{\frac{n}{2}}$ vertices. 
At least one of these stars does not contain $v_1$ as a vertex, and so none of the edges of this star are red.
Let the center of this star be $v_2$, and color all edges of this star red. 

In general, in step $i$ we consider the forest $F_i$ and let $F_i'$ be a subgraph of $G$ isomorphic to $F_i$.  The forest $F_i'$ consists of $i$ vertex-disjoint stars, each composed of $\floor{\frac{n}{i}}$ vertices. 
At least one of these stars does not contain a vertex from the set $\{v_1, v_2, \dots, v_{i-1}\}$; therefore, its edges are still not colored. 
We let the center of one such star be $v_i$ and color all of its edges red.

At the end of this procedure, we have 
\[
\sum_{k=1}^{\floor{\frac{n}{2}}}\left(\floor{\frac{n}{k}}-1 \right)>\sum_{k=1}^{\floor{\frac{n}{2}}}\left(\frac{n}{k}-2 \right)>
\int_{1}^{\floor{\frac{n}{2}}}\left(\frac{n}{k}-2 \right)\,dk
= n \log{n} -O(n)
\]
red edges. Thus we have $\s^*(n) \ge n\log{n}-O(n)$, completing the proof.

\section{Acknowledgements}
Gy\H{o}ri and Salia were supported by the National Research, Development and Innovation Office NKFIH, grant K132696.
Li was supported by NSFC (12071370) and Shaanxi Fundamental Science Research Project for Mathematics and Physics (22JSZ009).
Tompkins was supported by the National Research, Development and Innovation Office, NKFIH, grant K135800.

\bibliographystyle{abbrv}
\bibliography{references.bib}



\end{document}